\documentclass[11pt]{article}
\usepackage[utf8]{inputenc}
\usepackage[margin=1in]{geometry}
\usepackage{amsmath,amssymb,amsthm}
\usepackage{graphicx}
\usepackage{booktabs}
\usepackage{tikz}
\usetikzlibrary{patterns}
\usepackage{xcolor}
\usepackage[colorlinks=true,linkcolor=blue,citecolor=blue,urlcolor=blue]{hyperref}
\usepackage{authblk}
\usepackage{enumitem}

\newcommand{\amax}{\operatorname{max}}
\newtheorem{proposition}{Proposition}
\newtheorem{remark}{Remark}

\title{\bfseries The pre-commitment best-choice problem:\\ exact finite-$n$ formulas}
\author{Marcos Costa Santos Carreira}
\affil{A5X \quad \texttt{marcoscscarreira@gmail.com}}
\date{First arXiv version: 2018 \quad Revised: 2021 \quad This version: 2026}

\begin{document}
\maketitle

\begin{abstract}
In the full-information best-choice problem of Gilbert and Mosteller (1966), $n$ i.i.d.\
uniform draws are observed in sequence and the player, knowing $n$, stops at one with the
goal of catching the overall maximum. The optimal rule is adaptive---accept a draw only if
it is a running maximum and exceeds a round-dependent threshold---and its win probability
converges to $0.580164\ldots$ We study the restricted, non-adaptive \emph{pre-commitment}
class, in which the thresholds are fixed in advance and the player stops at the first draw
above its threshold, with no running-maximum check. Classifying each round as a win, a
false positive, a false negative, or a continuation, we derive exact finite-$n$ formulas
for all four probabilities---for a general non-increasing threshold vector, and in closed
form through the digamma function for a single repeated threshold---together with the
optimal thresholds and a constant-free approximation, all validated by simulation. We show
the class is strictly suboptimal for every $n\ge2$; a heuristic scaling-limit analysis
places its asymptotic win probability near $0.562$ (the single-threshold case gives
$0.51735$), below the optimal $0.580164$.
\end{abstract}

\section{Introduction}

The class of problems known as ``secretary problems'' has been studied extensively in
optimal stopping theory. The basic structure involves observing $n$ elements
sequentially, one at a time, and selecting one with an irrevocable decision made
immediately after each observation. In the \emph{full-information} version, the
observer learns the exact value of each draw and knows the underlying distribution,
in contrast to the \emph{no-information} version, in which only relative ranks are
revealed. The history and many variants are surveyed by Ferguson~\cite{Ferguson1989}.

In ``Recognizing the Maximum of a Sequence,'' Gilbert and Mosteller~\cite{GM1966}
analyse the full-information game in which $n$ measurements drawn from a uniform
distribution on $[0,1]$ are observed in sequence and the player must stop at the
overall maximum. They show that the optimal policy is an \emph{adaptive threshold}
rule. Writing $i = n-r$ for the number of rounds remaining at round $r$, the policy
accepts the current draw if and only if
\begin{enumerate}[label=(\roman*),nosep]
  \item the draw is a \emph{candidate}, i.e.\ a running maximum
        $a_r > \amax(a_1,\dots,a_{r-1})$, \emph{and}
  \item its value exceeds an indifference number $d_i$ that depends only on $i$.
\end{enumerate}
The indifference numbers solve
\begin{equation}\label{eq:indiff}
  \sum_{j=1}^{i}\binom{i}{j}\frac{d^{\,i-j}(1-d)^{j}}{j}=d^{\,i},
\end{equation}
giving $d_1=0$, $d_2=\tfrac12$, $d_3=0.68990,\dots$ The reason the candidate
restriction is correct is that a draw which is \emph{not} a running maximum is already
known not to be the overall maximum, so stopping at it can never win; the optimal rule
therefore never does. Under this rule, the overall win probability decreases in $n$
and converges to a constant for which Samuels~\cite{Samuels1982,Samuels1991} and
Gnedin and Miretskiy~\cite{GnedinMiretskiy2007} give the explicit asymptotic value
\begin{equation}\label{eq:opt}
  P_{\text{opt}} \;\xrightarrow[n\to\infty]{}\; 0.580164\ldots
\end{equation}

\paragraph{The pre-commitment variant.}
Gilbert and Mosteller also consider a simpler, \emph{non-adaptive} family of rules in
which the player fixes the thresholds in advance and stops at the first draw exceeding
its threshold---ignoring the running-maximum (candidate) check (ii). We call this the
\emph{pre-commitment} strategy. For a single threshold $k$ used in all rounds, Gilbert
and Mosteller derive the asymptotic win probability
\begin{equation}\label{eq:pcsingle}
  P_{\text{win}}^{k,\lim}=\sum_{i=1}^{\infty}\frac{e^{-\mu}\mu^{i}}{i!\,i},
  \qquad \mu\approx 1.503,\qquad P_{\text{win}}^{k,\lim}\approx 0.51735 .
\end{equation}
The pre-commitment strategy is interesting precisely because it is \emph{simple}: the
decision at round $r$ depends only on the current draw and a number fixed before the
game starts. It requires no memory of past draws and no online comparison, which makes
it attractive as a rule of thumb and as a building block for models with richer
payoffs. It is also, by construction, \emph{suboptimal}: because it can accept a draw
that is not a running maximum, it incurs avoidable losses. The purpose of this paper is
to give the exact finite-$n$ theory of this class and to quantify and explain how far
short of the optimum~\eqref{eq:opt} it falls. The explanation is geometric: a
pre-commitment rule can cut the unit hypercube only with axis-aligned slabs $a_r\ge k_r$,
never with the diagonal $a_r>\amax(a_1,\dots,a_{r-1})$ that defines a candidate, and that
missing cut is exactly the source of the gap.

\paragraph{A note on a previous version.} An earlier version of this work
(\,\href{https://arxiv.org/abs/1805.11556}{arXiv:1805.11556}\,) presented these exact
formulas but framed them as a \emph{correction} of an alleged error in Gilbert and
Mosteller, on the grounds that conditioning on having continued raises the expected
maximum of the remaining draws. That framing was mistaken. Future draws are
independent of past draws, so the probability that the next $i$ draws all fall below a
value $x$ is exactly $x^{i}$, as Gilbert and Mosteller assume; the elevated
\emph{conditional} maximum observed in simulation is a selection effect of conditioning
on ``the maximum has not yet appeared,'' which the optimal candidate rule already
accounts for. The decisive check is numerical: the same thresholds $d_i$
that give $P_{\text{win}}\approx 0.5837$ under the adaptive candidate rule give only
$\approx 0.5403$ when applied under the pre-commitment rule (Section~\ref{sec:sim},
Table~\ref{tab:n100}). The formulas derived here are correct---but they solve the
pre-commitment problem, not the optimal one. We therefore present them, in this
version, as what they are.

\paragraph{Contribution.} For the pre-commitment class we provide:
exact closed-form expressions for the per-round True Positive, False Positive, False
Negative and continuation probabilities, for general threshold vectors
(Section~\ref{sec:general}) and for the identical-threshold case
(Section~\ref{sec:identical}); the optimal pre-commitment thresholds and a
constant-free approximation (Section~\ref{sec:optimal}); a geometric account of the
gap to the optimum in terms of axis-aligned versus diagonal partitions of the unit
hypercube (Section~\ref{sec:geometry}); and simulation validation
(Section~\ref{sec:sim}).

\section{The pre-commitment strategy and four outcomes}\label{sec:setup}

Let $a_1,\dots,a_n$ be i.i.d.\ uniform on $[0,1]$. Fix a non-increasing threshold
vector $k=(k_1,\dots,k_{n-1},k_n)$ with $k_n=0$ (the last draw is always accepted). The
\emph{pre-commitment} rule stops at the first round $r$ with $a_r\ge k_r$. Because the
acceptance test ignores whether $a_r$ is a running maximum, each round yields exactly
one of four outcomes:
\begin{itemize}[nosep]
  \item \textbf{True Positive (Win):} accept and the accepted draw is the overall
        maximum, $a_r\ge k_r \wedge a_r=\amax(a_1,\dots,a_n)$;
  \item \textbf{False Positive (Loss):} accept but the draw is not the overall maximum,
        $a_r\ge k_r \wedge a_r\neq \amax(\cdot)$;
  \item \textbf{False Negative (Loss):} do not accept although the draw \emph{is} the
        overall maximum, $a_r< k_r \wedge a_r=\amax(\cdot)$;
  \item \textbf{Continuation:} do not accept and the draw is not the overall maximum,
        $a_r< k_r \wedge a_r\neq \amax(\cdot)$.
\end{itemize}
The False Positive category is exactly what the optimal candidate rule avoids: it
arises only because pre-commitment accepts draws above threshold that a single
comparison with the running maximum would have rejected. We write $PW$, $PFP$, $PFN$
and $PC$ for the four per-round probabilities and seek closed forms in $n$, $r$ and
$k$.

\section{Geometric visualisation}\label{sec:geometry}

For $n=2$ the pair $(a_1,a_2)$ fills the unit square. The diagonal from $(0,0)$ to
$(1,1)$ separates $a_1>a_2$ (below) from $a_1<a_2$ (above); the vertical line
$a_1=k_1$ separates accept ($a_1\ge k_1$) from continue. The four regions of
Section~\ref{sec:setup} are shown in Figure~\ref{fig:square}.

\begin{figure}[h]
\centering
\begin{tikzpicture}[scale=4.2]
  \fill[red!55]    (0.6,0) -- (1,0) -- (1,1) -- (0.6,0.6) -- cycle;
  \fill[brown!55]  (0.6,0.6) -- (1,1) -- (0.6,1) -- cycle;
  \fill[blue!45]   (0,0) -- (0.6,0) -- (0.6,0.6) -- cycle;
  \fill[yellow!65] (0,0) -- (0.6,0.6) -- (0.6,1) -- (0,1) -- cycle;
  \draw[thick] (0,0) rectangle (1,1);
  \draw[thick] (0,0) -- (1,1);
  \draw[thick,blue!60!black,dashed] (0.6,0) -- (0.6,1);
  \node at (0.85,0.30) {\small\textbf{Win}};
  \node at (0.78,0.92) {\small FP};
  \node at (0.30,0.12) {\small FN};
  \node at (0.22,0.72) {\small Continue};
  \node[below] at (0.6,-0.02) {\small $k_1$};
  \node[below] at (1,-0.02) {\small $1$};
  \node[left] at (0,1) {\small $a_2$};
  \node[below] at (1.0,1.04) {};
  \draw[->] (0,0)--(1.08,0) node[right]{\small $a_1$};
  \draw[->] (0,0)--(0,1.08);
\end{tikzpicture}
\caption{The four outcomes for $n=2$ with threshold $k_1$. The accept boundary
$a_1=k_1$ is a vertical (axis-aligned) line; the win/loss boundary $a_1=a_2$ is the
diagonal. Pre-commitment can place the vertical cut anywhere but cannot bend it onto
the diagonal.}
\label{fig:square}
\end{figure}

The picture generalises. For $n=3$ the draws fill the unit cube and the optimal rule's
candidate condition $a_r>\amax(a_1,\dots,a_{r-1})$ carves the cube along its diagonal,
producing non-rectangular regions; in particular, once the round-1 False Negative
region is removed, the conditional distribution entering round 2 is no longer the
uniform square of the $n=2$ problem. The pre-commitment rule, by contrast, only ever
inserts \emph{axis-aligned slabs} $a_r\ge k_r$. These tile neatly, which is exactly why
closed forms exist; but they cannot reproduce the diagonal cut, which is the structural
source of the suboptimality quantified below. The False Positive mass in
Figure~\ref{fig:square} (the brown triangle) is precisely the probability the optimal
rule reclaims by replacing the vertical cut with the diagonal one.

\section{Exact formulas for general thresholds}\label{sec:general}

Throughout this section the threshold vector is assumed \emph{non-increasing},
$k_1\ge k_2\ge\cdots\ge k_{n-1}\ge k_n=0$; this is the only regime in which the
integrated regions tile in the form below, and it entails no loss of generality for the
optimisation, since the optimal pre-commitment thresholds are themselves non-increasing.
Integrating the four regions over the unit hypercube (automated symbolically for
$n\le 6$) yields a stable pattern. Set $PC(n,0)=1$. The False Negative probability
depends only on the threshold of its own round,
\begin{equation}\label{eq:pfn}
  PFN(n,r)=\frac{k_r^{\,n}}{n},\qquad 1\le r\le n .
\end{equation}
The win and continuation probabilities are tied to the continuation probability of the
previous round. The justification is \emph{exchangeability}: conditioned on the event
that the game has continued through round $r-1$, the $n-r+1$ still-unseen draws occupy
positions $r,\dots,n$ that are exchangeable, so the overall maximum of those positions
is equally likely to sit at any one of them. The chance that it sits at the current
position $r$ is therefore $1/(n-r+1)$, and the win-plus-false-negative mass available at
round $r$ is the previous continuation mass spread over those $n-r+1$ positions,
\begin{equation}\label{eq:pw_rec}
  PW(n,r)=\frac{PC(n,r-1)}{\,n-r+1\,}-PFN(n,r),\qquad 1\le r\le n,
\end{equation}
and the False Positive probability is the residual
\begin{equation}\label{eq:pfp}
  PFP(n,r)=\bigl(PC(n,r-1)-(PW(n,r)+PFN(n,r))\bigr)-PC(n,r).
\end{equation}
Everything therefore reduces to a closed form for the continuation probability
$PC(n,r)$. Writing the integrated tables for $n$ up to $6$ in a common form reveals the
decomposition
\begin{equation}\label{eq:pc}
  PC(n,r)=\prod_{j=1}^{r}k_j-\sum_{i=1}^{r}\Bigl(g(i,n,r)\prod_{j=1}^{r}k_j^{\,f(i,j,n,r)}\Bigr),
\end{equation}
with exponent and coefficient arrays
\begin{equation}\label{eq:fg}
  f(i,j,n,r)=\begin{cases} n-r+i & i=j,\\[2pt] 1 & i<j,\\[2pt] 0 & i>j,\end{cases}
  \qquad
  g(i,n,r)=\frac{n-r}{(n-r+i-1)(n-r+i)} ,
\end{equation}
and $\sum_{j=1}^{r}f(i,j,n,r)=n$ for every $i,r$. Substituting~\eqref{eq:pc}
into~\eqref{eq:pw_rec} gives the per-round win probability for a non-increasing threshold
vector,
\begin{equation}\label{eq:pw}
  PW(n,r)=\frac{1}{n-r+1}\Biggl(\prod_{j=1}^{r-1}k_j
        -\sum_{i=1}^{r-1}g(i,n,r-1)\prod_{j=1}^{r-1}k_j^{\,f(i,j,n,r-1)}\Biggr)-\frac{k_r^{\,n}}{n},
\end{equation}
and the overall win probability is
\begin{equation}\label{eq:pwn}
  PW(n)=\sum_{r=1}^{n}PW(n,r).
\end{equation}
Equation~\eqref{eq:pwn} is the objective maximised over $k_1,\dots,k_{n-1}$ to obtain
the best pre-commitment strategy. We are not aware of a prior closed-form treatment of
this non-adaptive multi-threshold class; the literature in
Section~\ref{sec:related} treats either the optimal adaptive rule or the single-threshold
asymptotics only.

\section{The identical-threshold case}\label{sec:identical}

Setting $k_1=\dots=k_{n-1}=k$ (and $k_n=0$) collapses~\eqref{eq:pw} to
\begin{equation}\label{eq:pwk}
  PW_k(n,r)=\frac{k^{\,r-1}}{n-r+1}-\frac{k^{n}}{n-r+1}
            +\begin{cases}0 & 1\le r\le n-1,\\[2pt] k^{n}/n & r=n,\end{cases}
\end{equation}
or, collecting the common factor $k^{n}$ (note $k^{n}/[(n-r+1)k^{\,n-r+1}]=k^{\,r-1}/(n-r+1)$),
\begin{equation}\label{eq:pwk2}
  PW_k(n,r)=k^{n}\!\left(\frac{1}{(n-r+1)\,k^{\,n-r+1}}-\frac{1}{n-r+1}\right)
            +\begin{cases}0 & 1\le r\le n-1,\\[2pt] k^{n}/n & r=n.\end{cases}
\end{equation}
The per-round values are non-negative, and summing over $r$ the overall win probability
admits the special-function form
\begin{equation}\label{eq:pwk_total}
  PW_k(n)=k^{n}\!\left(\sum_{r=1}^{n}\frac{1}{(n-r+1)\,k^{\,n-r+1}}-\bigl(\psi^{(0)}(n)+\gamma\bigr)\right),
\end{equation}
where $\psi^{(0)}$ is the digamma function and $\gamma$ is the Euler--Mascheroni
constant (so $\psi^{(0)}(n)+\gamma=\sum_{m=1}^{n-1}1/m$, the harmonic number). Optimising
$PW_k(n)$ over $k$ is fast; the optimum sits near $k^\star\approx 1-1.5/n$ for large
$n$, and $PW_k(n)\to 0.51735$, recovering the Gilbert--Mosteller single-threshold value
\eqref{eq:pcsingle}.

\section{Optimal thresholds and a constant-free approximation}\label{sec:optimal}

\begin{figure}[h]
\centering
\includegraphics[width=0.78\textwidth]{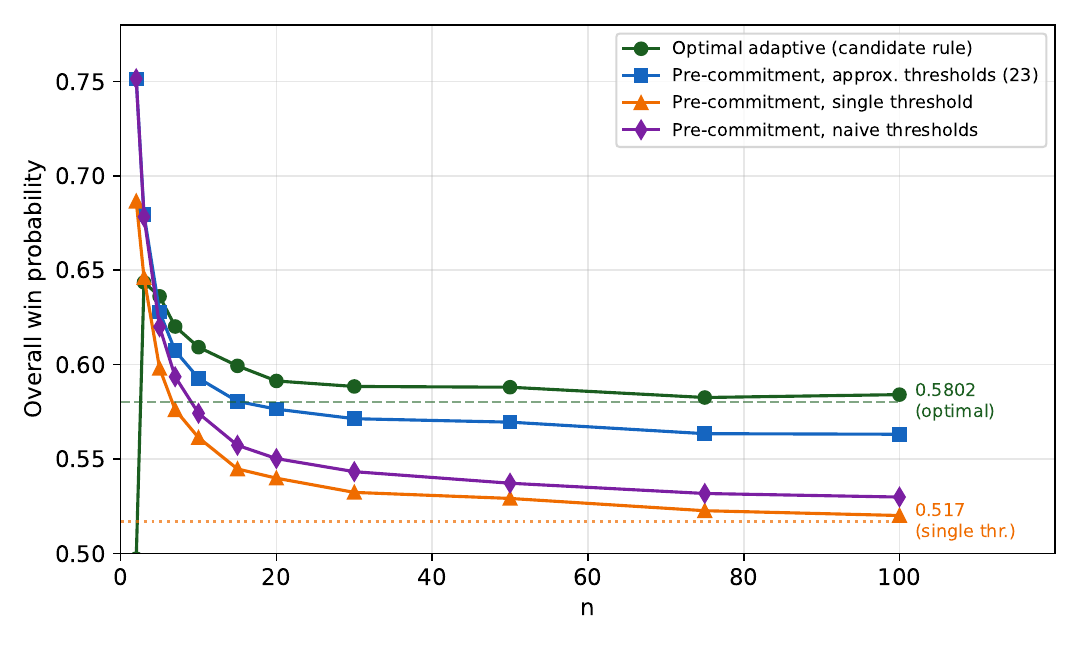}
\caption{Overall win probability versus $n$. The optimal adaptive (candidate) rule
converges to $0.580164$; the best pre-commitment rule (approximate thresholds
\eqref{eq:kapp}) converges to $\approx0.56$; the single identical threshold converges to
$0.517$. Markers are simulation estimates ($1.2\times10^{5}$ runs per point).}
\label{fig:winprob}
\end{figure}

\begin{figure}[h]
\centering
\includegraphics[width=0.78\textwidth]{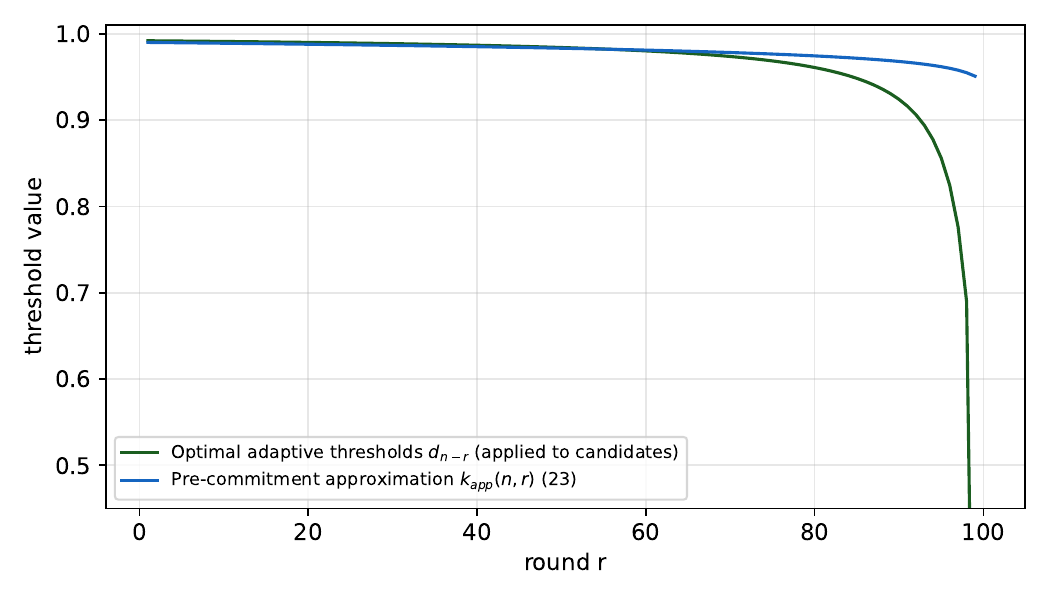}
\caption{Threshold profiles for $n=100$: the optimal adaptive indifference numbers
$d_{n-r}$ (applied only to candidates) versus the pre-commitment approximation
$k_{\mathrm{app}}(n,r)$ of \eqref{eq:kapp}.}
\label{fig:thresholds}
\end{figure}

\begin{figure}[h]
\centering
\includegraphics[width=0.78\textwidth]{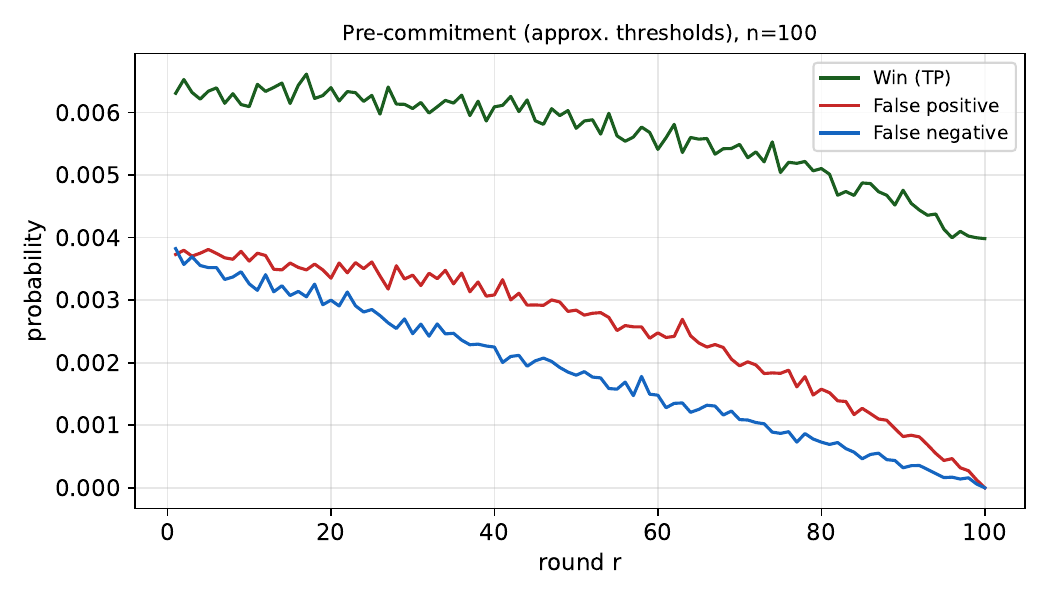}
\caption{Per-round outcome probabilities for the approximate pre-commitment rule,
$n=100$: slowly decaying wins, False Positives decaying from the start, and roughly
linearly decaying False Negatives.}
\label{fig:decomp}
\end{figure}

Maximising~\eqref{eq:pwn} over the full threshold vector gives the optimal
pre-commitment strategy. The optimal $k_r$ increase with $n$ and, for fixed $n$,
decrease across rounds; the per-round win probability is almost flat early and decays
slowly, while False Positives decay from the start and False Negatives decay roughly
linearly toward zero. Empirically the optimal thresholds are close to linear in
$\log(n-r)$, which suggests the constant-free approximation
\begin{equation}\label{eq:kapp}
  k_{\mathrm{app}}(n,r)=\Bigl(1-\tfrac1n\Bigr)+\tfrac1n\,\log\!\Bigl(\tfrac{n-r}{n}\Bigr).
\end{equation}
The first term is the naive choice $1-1/n$; the $1/n$ multiplying the logarithm comes
from the observed decay of the slopes with $n$. The choice $1-1/n$ on the first round
also has the property that the round-1 False Positive and False Negative probabilities
are equal. Despite having no fitted constants, \eqref{eq:kapp} performs nearly as well
as full numerical optimisation: for $n=100$ it gives $PW\approx 0.5634$.
Figure~\ref{fig:winprob} compares the overall win probability of the optimal adaptive
rule, the optimal/approximate pre-commitment rule, the single-threshold pre-commitment
rule and the naive rule as functions of $n$; Figure~\ref{fig:thresholds} shows the
threshold profiles for $n=100$; and Figure~\ref{fig:decomp} shows the per-round
decomposition for the approximate pre-commitment rule.

\section{Simulation}\label{sec:sim}

We generate $10^{6}$ runs of $100$ draws and take the first $n$ of each run. For each
run and each rule we record whether the outcome is a Win, a False Positive or a False
Negative, classified by round, and compare the empirical frequencies
with~\eqref{eq:pw}, \eqref{eq:pfp} and \eqref{eq:pfn}. The predictions match the
simulations to four decimal places for all tested $n$, confirming that the formulas
are exact for the pre-commitment problem.

The same simulations make the central distinction concrete. We apply
\emph{one and the same} set of Gilbert--Mosteller thresholds in two ways: under the
optimal \emph{adaptive} rule (accept a draw only if it is a candidate and exceeds the
threshold) and under the \emph{pre-commitment} rule (accept the first draw above the
threshold). Table~\ref{tab:n100} reports $n=100$.

\begin{table}[h]
\centering
\begin{tabular}{lcc}
\toprule
Rule applied with G\&M thresholds $d_{n-r}$ & Win probability & Note\\
\midrule
Adaptive (candidate check, the optimal rule) & $0.5837$ & $\to 0.580164$ \\
Pre-commitment (no candidate check)          & $0.5403$ & this paper's rule \\
\midrule
Pre-commitment, optimal thresholds           & $0.5651$ & \\
Pre-commitment, approximation \eqref{eq:kapp}& $0.5634$ & \\
Pre-commitment, single threshold             & $0.5218$ & $\to 0.51735$ \\
\bottomrule
\end{tabular}
\caption{Simulated win probabilities, $n=100$. The gap between the first two rows
($0.5837$ vs.\ $0.5403$) is entirely the candidate check: the same thresholds, applied
adaptively, win far more often. The best pre-commitment rule ($0.5651$) still falls
short of the optimal $0.5837$.}
\label{tab:n100}
\end{table}

The best non-adaptive pre-commitment strategy thus recovers roughly three quarters of
the gap between the single-threshold rule and the optimum, but it cannot close it: the
remaining shortfall is exactly the False Positive mass that only the diagonal candidate
cut of Section~\ref{sec:geometry} can remove.

\section{Suboptimality and a scaling-limit heuristic}\label{sec:subopt}

So far the gap between the best pre-commitment rule and the optimum has been shown
numerically (Section~\ref{sec:sim}) and explained geometrically
(Section~\ref{sec:geometry}). We now record it as a theorem. Let $v(n)$ be the optimal
win probability over all stopping times adapted to the natural filtration of
$a_1,\dots,a_n$, and let $v_{pc}(n)=\sup_k PW(n)$ be the supremum over non-increasing
pre-commitment threshold vectors. It is classical that $v(n)$ decreases to
$v_\infty=0.580164\ldots$

\begin{proposition}\label{prop:subopt}
For every $n\ge 1$, $v_{pc}(n)\le v(n)$, and for every $n\ge 2$ the inequality is strict.
Consequently $\limsup_{n\to\infty} v_{pc}(n)\le v_\infty=0.580164\ldots$
\end{proposition}

\begin{proof}
\textit{(Subclass bound and a coupling.)} A pre-commitment rule ``stop at the first $r$
with $a_r\ge k_r$'' is a stopping time, so the class is a subset of all stopping rules
and $v_{pc}(n)\le v(n)$. To see where mass is lost, fix a threshold vector $k$ and couple
the pre-commitment rule $\sigma$ with the rule $\tilde\sigma$ that uses the \emph{same}
thresholds but stops only at a draw that is simultaneously above its threshold and a
candidate (running maximum), taking the last draw if no such round occurs. The two agree
until the first crossing $\tau=\min\{r:a_r\ge k_r\}$. If $a_\tau$ is a candidate, both
stop at $\tau$ with identical outcome. If $a_\tau$ is \emph{not} a candidate then
$a_\tau<\max(a_1,\dots,a_{\tau-1})$, so $a_\tau$ is not the overall maximum, $\sigma$
loses with certainty, and $\tilde\sigma$ continues. Hence, with
$F=\{a_\tau\text{ not a candidate}\}$,
\begin{equation}\label{eq:coupling}
  W(\tilde\sigma)-W(\sigma)=\Pr(F)\,\Pr(\tilde\sigma\text{ wins}\mid F)\ \ge\ 0 .
\end{equation}
Taking the supremum over $k$ re-gives $v_{pc}(n)\le v(n)$, and $\limsup_n v_{pc}(n)\le
\lim_n v(n)=v_\infty$.

\textit{(Strictness for finite $n$.)} The optimum $v(n)$ is attained only by rules that
never stop at a non-candidate, since stopping at a non-candidate wins with probability
$0$ and is strictly dominated by continuing (future draws are independent of the past and
carry positive winning probability). Every non-increasing pre-commitment rule with
$k_1<1$ stops at a non-candidate with positive probability: at $n=2$ the region
$\{k_1\le a_1,\ a_1<a_2\}$ has positive area, and for general $n$ the region
$\{k_r\le a_r<\max(a_1,\dots,a_{r-1})\}$ at any round $r\ge 2$ has positive measure
whenever $k_r<1$. (If $k_1=1$, the rule cannot stop early and is suboptimal directly.)
Thus no pre-commitment rule attains $v(n)$ and $v_{pc}(n)<v(n)$ for all $n\ge 2$.
\end{proof}

The remaining question---whether the limit $v_{pc}^\infty:=\lim_n v_{pc}(n)$ is itself
\emph{strictly} below $v_\infty$ or merely approaches it---is best read through the
standard scaling limit. We treat this limit \emph{heuristically}: the derivation that
follows is formal and the resulting constant is obtained by numerical optimisation. The
finite-$n$ conclusions of Proposition~\ref{prop:subopt} are rigorous and do not depend on
any of it.

\subsection{A scaling-limit heuristic}\label{sec:scaling}

Rescale each draw by its distance below the top, $y_r:=n(1-a_r)$, and time by $t:=r/n$.
The point process $\Pi_n=\{(r/n,\,y_r):1\le r\le n\}$ on $[0,1]\times[0,\infty)$ converges
in distribution, as $n\to\infty$, to a Poisson point process $\Pi$ with intensity
$\mathrm{d}t\,\mathrm{d}y$ (Lebesgue): each $a_r$ is uniform, so $\Pr(y_r\in\mathrm{d}y)=
\mathrm{d}y/n$ on $[0,n]$, and the $n$ independent points over $t\in[0,1]$ give the
unit-rate limit. A draw is \emph{above threshold}, $a_r\ge k_r$, exactly when its point
lies below the curve $y\le n(1-k_r)$; we therefore scale thresholds as
\begin{equation}\label{eq:betadef}
  n\bigl(1-k_{\lfloor tn\rfloor}\bigr)\longrightarrow \beta(t),
\end{equation}
for a non-decreasing profile $\beta:[0,1)\to[0,\infty)$ (non-increasing $k$ corresponds
to non-decreasing $\beta$). Under \eqref{eq:betadef} the per-round false-negative mass
satisfies $k_r^{\,n}\to e^{-\beta(t)}$, recovering $\mathrm{FN}_{\mathrm{tot}}\to
\int_0^1 e^{-\beta(t)}\,\mathrm{d}t$.

In the limit the overall maximum is the point of $\Pi$ with smallest $y$; its height
$Y$ is $\mathrm{Exp}(1)$ (the number of points with $y<\varepsilon$ is Poisson with mean
$\varepsilon$, so $\Pr(Y>\varepsilon)=e^{-\varepsilon}$) and its time $T$ is independent
and uniform on $[0,1]$. The pre-commitment rule stops at the first $t$ with a point below
$\beta$; it wins iff that first crossing is the global minimum $Y$, i.e.\ iff $Y\le
\beta(T)$ \emph{and} no point with time $<T$ lies below the curve. Conditioning on the
minimal point at $(T,Y)$, the remaining points form an independent Poisson process on
$\{y>Y\}$, and the mean number of earlier crossings is $\int_0^T(\beta(s)-Y)^+\mathrm{d}s$.
Hence, formally, the limiting win probability is the functional
\begin{equation}\label{eq:Vfunc}
  V(\beta)=\int_0^1\!\!\int_0^{\beta(t)} e^{-y}\,
            \exp\!\Bigl(-\!\int_0^t(\beta(s)-y)^+\,\mathrm{d}s\Bigr)\,\mathrm{d}y\,\mathrm{d}t ,
\end{equation}
and $v_{pc}^\infty=\sup_{\beta\,\uparrow} V(\beta)$ over non-decreasing profiles.

\begin{remark}[consistency check]
For a constant profile $\beta(t)\equiv\mu$ (a single threshold $k=1-\mu/n$),
\eqref{eq:Vfunc} reduces to $V(\mu)=\int_0^1\!\int_0^{\mu}e^{-y}e^{-(\mu-y)t}\,
\mathrm{d}y\,\mathrm{d}t$, whose maximiser is $\mu^\ast=1.5029$ with
$V(\mu^\ast)=0.517351$. This is exactly Gilbert and Mosteller's single-threshold
constant~\eqref{eq:pcsingle}, $\sum_{i\ge1}e^{-\mu}\mu^i/(i!\,i)$, confirming the
functional.
\end{remark}

Two features of $V$ are worth recording, one rigorous and one numerical.

\emph{(i) Per-profile strictness.} The candidate-completion coupling \eqref{eq:coupling}
carries over to the limit: for a fixed profile $\beta$, the adaptive rule that accepts only
\emph{records} (running minima of $y$) lying below $\beta$ is admissible and stops no worse
than the pre-commitment rule, so $V(\beta)\le v_\infty$; and the inequality is strict,
$V(\beta)<v_\infty$, for every profile that crosses with positive probability, because the
pre-commitment rule is forced to accept non-record crossings (sure losses) with positive
probability. Thus no \emph{single} profile attains the adaptive value.

\emph{(ii) Stationarity and the numerical optimum.} A formal first variation of
\eqref{eq:Vfunc}, perturbing $\beta\mapsto\beta+\varepsilon\eta$ (the upper limit
contributing $e^{-\beta(s)}$ since $J(s,\beta(s))=0$, the interior contributing the double
integral below), yields the stationarity condition
\begin{equation}\label{eq:EL}
  e^{-\beta^\star(s)}=\int_s^1\!\!\int_0^{\beta^\star(s)}
       e^{-y}\,e^{-J(t,y)}\,\mathrm{d}y\,\mathrm{d}t ,
  \qquad J(t,y)=\int_0^t(\beta^\star(u)-y)^+\,\mathrm{d}u ,
\end{equation}
balancing, at each time, the gain from lowering the threshold (catching a would-be false
negative) against the premature crossings it admits. Optimising \eqref{eq:Vfunc}
numerically---equivalently, extrapolating the exact finite-$n$ optima---gives
\[
   v_{pc}^\infty\approx 0.562,\qquad \beta^\star(t)\approx 0.90-0.94\log(1-t),
\]
comfortably below $v_\infty=0.580164\ldots$, while the single-threshold restriction returns
Gilbert and Mosteller's $0.51735$ (the series $\sum_{i\ge1}e^{-\mu}\mu^i/(i!\,i)$ at
$\mu=1.503$).

\begin{remark}
We deliberately do \emph{not} present this scaling limit as a theorem. A rigorous
treatment would have to justify the convergence to the Poisson functional, the
compactness of the admissible profile class, and the attainment of the supremum---enough
analysis to form a separate paper, which we leave to future work. Nor do we expect an
elementary closed form for $v_{pc}^\infty$: even the adaptive optimum
$v_\infty=0.580164\ldots$ has none, being defined through its own limiting integral
equation~\cite{GnedinMiretskiy2007}. What the heuristic does add is a clean deterministic
account of \emph{why} the gap persists. The term $\int_0^t(\beta(s)-y)^+\mathrm{d}s$ in
\eqref{eq:Vfunc} is the expected number of \emph{premature} crossings---points below the
fixed curve that need not be running minima---which is exactly the false-positive mass the
diagonal candidate test $a_r>\max(a_1,\dots,a_{r-1})$ would suppress and a fixed profile
cannot. The finite-$n$ evidence (Table~\ref{tab:n100}, Figure~\ref{fig:winprob}) and this
limit agree on $v_{pc}^\infty\approx 0.56<v_\infty$.
\end{remark}

\section{Related work}\label{sec:related}

The optimal full-information best-choice problem is settled. Gilbert and
Mosteller~\cite{GM1966} pose it and solve the optimal adaptive rule; Samuels gives exact
finite-$n$ solutions~\cite{Samuels1982} and the survey treatment~\cite{Samuels1991};
Gnedin and Miretskiy~\cite{GnedinMiretskiy2007} derive the explicit asymptotic winning
rate $0.580164$, and Gnedin~\cite{Gnedin1996} establishes structural properties of the
optimal threshold sequence. Ferguson~\cite{Ferguson1989} surveys the field, and
Bruss's odds theorem~\cite{Bruss2000} provides a general last-success stopping tool
that underlies the per-round bookkeeping used here.

The closest related work on non-adaptive thresholds comes not from classical secretary
theory but from the prophet and prophet-secretary literature, where static or blind
threshold policies are studied for competitive guarantees against an offline prophet.
Correa, Saona and Ziliotto~\cite{BlindStrategies} analyse non-increasing \emph{blind}
threshold sequences that depend only on the distribution of the maximum; Arnosti and
Ma~\cite{StaticThreshold} prove tight guarantees for static (fixed) threshold policies;
and Goldenshluger, Malinovsky and Zeevi~\cite{SingleThreshold} characterise optimal
single-threshold rules through a zero-sum game on the unit square. These models differ
from the present objective: they maximise expected accepted value relative to the prophet
benchmark, typically under heterogeneous distributions, whereas here the objective is the
exact probability of selecting the sample maximum in the uniform full-information
best-choice problem. For that objective the non-adaptive class has previously been
treated only in the single-threshold asymptotic regime~\eqref{eq:pcsingle}. Existing
static-threshold results are thus conceptually adjacent but do not appear to give the
finite-$n$ TP/FP/FN/continuation formulas derived here for the multi-threshold
pre-commitment class. To the author's knowledge, exact finite-$n$ formulas for the
general non-increasing multi-threshold pre-commitment rule, together with this
four-outcome decomposition, do not appear to have been given previously.

\section{Conclusion}

We have given an exact finite-$n$ theory of the non-adaptive pre-commitment strategy
for the full-information best-choice problem: closed forms for the four per-round
outcome probabilities, a digamma expression for the identical-threshold case, the
optimal thresholds, and a simple constant-free approximation that nearly matches them.
The strategy is valuable for its simplicity---no memory, no online comparison---and the
explicit four-outcome decomposition extends readily to problems with round-dependent
payoffs or asymmetric costs for the two kinds of loss.

The strategy is also, unavoidably, suboptimal: Proposition~\ref{prop:subopt} shows the
best pre-commitment win probability is strictly below the optimum for every $n\ge2$, and a
heuristic scaling-limit analysis (Section~\ref{sec:scaling}) locates the limit near
$0.562$, below $0.580164$, with the single-threshold sub-case at $0.51735$. The reason is
structural and geometric:
pre-commitment partitions the
unit hypercube with axis-aligned slabs and cannot represent the diagonal cut
$a_r>\amax(a_1,\dots,a_{r-1})$ that defines a candidate. That diagonal is exactly what
the optimal adaptive rule of Gilbert and Mosteller exploits, and applying their
thresholds with the candidate check restores the optimal value. Building the
geometric picture before the derivation---rather than after---would have made the
distinction between the two rules visible from the outset.

\paragraph{Acknowledgements.} This problem was posted on \emph{Cut the Knot} by the late
Alexander Bogomolny, who passed away in 2018. The author is deeply grateful to him for
creating a place to try, fail, and try again, and dedicates this work to his memory; his
generosity and his love of a good problem are missed. The author is also grateful to
Steven Strogatz, Michael Wiener and Hank Tijms for their comments and inputs, and to
Alexander Gnedin---whose work
\cite{Gnedin1996,GnedinMiretskiy2007} underpins the asymptotics of
Section~\ref{sec:subopt}---for his inputs. Code for the integration, figures and
simulations is available at
\href{https://github.com/MarcosCarreira/AnExactSolution}{github.com/MarcosCarreira/AnExactSolution}.

\paragraph{Use of AI tools.} The author used Claude (Anthropic) as an assistant in
preparing this revision: reframing the exposition, surveying and positioning the related
literature, generating figures, performing numerical and symbolic verification of the
formulas, and drafting the scaling-limit argument of Section~\ref{sec:subopt}. All
results were checked and are the responsibility of the author.

\end{document}